\documentclass[a4 paper, 12pt]{article}

\usepackage{amsfonts}
\usepackage{amssymb}
\usepackage{amsmath}
\usepackage{fancyhdr}
\usepackage{eufrak}
\usepackage{latexsym}
\usepackage{amsbsy}
\usepackage[all]{xy}
\usepackage[latin1]{inputenc}
\usepackage[T1]{fontenc}        
\usepackage{amsmath}

\newenvironment{theorem}[1]
{\vskip 2mm\noindent \textbf{Theorem #1}  \it}{\vskip 2mm}

\newenvironment{proposition}[1]
{\vskip 2mm\noindent \textbf{Proposition #1}  \it}{\vskip 2mm}

\newtheorem{defn}{Definition}[section]
\newtheorem{thm}[defn]{Theorem}
\newtheorem{prop}[defn]{Proposition}
\newtheorem{lem}[defn]{Lemma}

\newenvironment{corollary}[1]
{\vskip 2mm\noindent \textbf{Corollary #1}  \it}{\vskip 2mm}

\newenvironment{claim}[1]
{\vskip 2mm\noindent \textbf{Claim #1}  \it}{\vskip 2mm}

\newenvironment{fact}[1]
{\vskip 2mm\noindent \textbf{Growth lemma #1}  \it}{\vskip 2mm}

\newenvironment{conj}[1]
{\vskip 2mm\noindent \textbf{Conjecture #1}  \it}{\vskip 2mm}

\newcommand{\proof}{\vskip 2mm \noindent {\textsc{Proof: }}\rm}
\newcommand{\proo}{\vskip 2mm \noindent {\textsc{Proof }}\rm}
\newcommand{\fin}{\hfill{\Large$\Box$}\\}
\newcommand{\finsec}{\hfill{\Large$\Box$}}

\newcommand{\al}{\alpha}

\newcommand{\ga}{\gamma}
\newcommand{\pp}{\nu \textrm{-a.e.}}
\newcommand{\si}{\sigma}
\newcommand{\Om}{\Omega}
\newcommand{\om}{\omega}
\newcommand{\epsi}{\epsilon}

\newcommand{\Sg}{\mathbb {S}}
\newcommand{\C}{\mathbb {C}}
\newcommand{\D}{\mathbb {D}}
\newcommand{\R}{\mathbb {R}}
\newcommand{\N}{\mathbb {N}}
\newcommand{\Z}{\mathbb {Z}}

\newcommand{\Pj}{\mathbb {C}\mathbb {P}}

\newcommand{\Ld}{\textsf {L}}

\newcommand{\cdb}{\textsf {Card }}
\newcommand{\orb}{\textsf {O}}

\newcommand{\Vol}{\textsf {Vol }}

\newcommand{\Ad}{\textsf {A}}
\newcommand{\pd}{\textsf {p}}

\newcommand{\hd}{\textsf {h}}
\newcommand{\vd}{\textsf {v}}
\newcommand{\wwd}{\textsf {w}}

\newcommand{\Id}{{\rm Id}}

\newcommand{\Lip}{{\rm Lip \, }}
\newcommand{\Jac}{{\rm Jac \,  }}
\newcommand{\dist}{{\textsf {dist}  }}

\def\su#1{\textrm {supp} (#1)}
\def\abs#1{\vert #1\vert}
\def\norm#1{\left\|\, #1\,\right\|}

\def\RRR{{\mathfrak R}}

\def\EE{{\cal E}}

\def\CC{{\cal C}}
\def\PP{{\cal P}}

\def\FF{{\cal F}}

\def\RR{{\cal R}}

\def\LL{{\cal L}}

\def\SS{{\cal S}}

\def\QQ{{\cal Q}}

\def\WW{{\cal W}}

\def\HH{{\cal H}}
\def\com{\ar@{}[rd]|{\circlearrowleft}}

\title {On the dimension of invariant measures of endomorphisms of $\Pj^k$}

\author{Christophe Dupont}
\date{ \today }

\begin{document}

\maketitle

\begin{abstract} Let $f$ be an endomorphism of $\Pj^k$ and $\nu$ be an
  $f$-invariant measure with positive Lyapunov exponents $(\lambda_1,\ldots,\lambda_k)$. We prove a
  lower bound for the pointwise
  dimension of $\nu$ in terms of the degree of $f$, the exponents  of
  $\nu$ and
  the entropy of
  $\nu$. In particular our result can be applied for the maximal entropy
  measure $\mu$. When
  $k=2$, it implies that the Hausdorff dimension of $\mu$ is estimated by $\dim_\HH \mu \geq {\log d \over  \lambda_1} + { \log d \over
  \lambda_2}$, which is half of the conjectured formula. Our method for proving these results consists in
  studying the distribution of the $\nu$-generic inverse
  branches of $f^n$  in $\Pj^k$. Our tools are a volume growth estimate for the
  bounded holomorphic polydiscs in $\Pj^k$ and a normalization theorem
  for the $\nu$-generic inverse branches of $f^n$.

\end{abstract}

\small{\noindent \emph{Key Words}: holomorphic dynamics, dimension theory.

}\normalsize

\section{Introduction}

Let $f$ be a smooth map acting on a compact Riemannian
manifold $M$ and $\nu$  be an $f$-invariant measure on $M$. By Young \cite{Y}, the pointwise
dimension of $\nu$ is defined by (provided the limit exists):
\begin{equation*}\label{exdi}
 \delta(x) = \lim_{r \to 0} \, {\log \nu (B_x(r)) \over  \log r},
\end{equation*}
where $B_x(r)$ is the ball in $M$ of center $x$ and radius $r$
(take $\liminf$ and $\limsup$ to define the lower and upper pointwise dimensions
$\underline \delta$ and $\bar \delta$). That function actually describes the geometrical behaviour of $\nu$ with
respect to the metric of $M$:  if $a \leq \underline
\delta \leq \bar \delta \leq b$ hold $\pp$, then the Hausdorff
dimension of $\nu$ also satisfies $a \leq \dim_\HH \nu \leq b$
\cite{Y}. Recall that $\dim_\HH \nu$ is defined as the infimum of the Hausdorff dimension
of the full $\nu$-measure borel subsets in $M$. In particular we have
$\dim_\HH \nu \leq \dim_\HH \su \nu$. We refer to the book of Pesin
\cite{P} for an introduction to dimension theory in dynamical
systems.   \\

Given a dynamical system $(M,f,\nu)$, one can expect relations between
the dimension of $\nu$, its Lyapunov exponents $\lambda_k \leq \ldots
\leq \lambda_1$ and its entropy $h_\nu$ (see \cite{Le}, \cite{P}). The situation has been completely
described when $f$ is a smooth diffeomorphism and $\nu$ is an
$f$-invariant hyperbolic measure (i.e. with no zero exponents).  Young
\cite{Y} first proved in the case of surfaces the formula $\delta =  h_\nu / \lambda_1 - h_\nu
 / \lambda_2$ $\pp$, where $\lambda_2 < 0 < \lambda_1$. In higher
 dimensions, Ledrappier-Young \cite{LY} established that the unstable
 pointwise dimension of $\nu$ satisfies $\pp$
\begin{equation}\label{ly} 
\delta^u = { h_1 \over  \lambda_1} + \sum_{i = 2}^u { h_i - h_{i-1}
  \over \lambda_i} ,  \end{equation}
 where  $h_1 \leq \ldots \leq h_u = h_\nu$ denote the conditional entropies of $\nu$ along
 the unstable manifolds $\WW^1 \subset \ldots \subset
 \WW^u$ (a similar formula holds for the stable dimension  $\delta^s$).
Later Barreira-Pesin-Schmeling \cite{BPS} proved the formula $\delta =
\delta^s + \delta^u$ $\pp$ by showing a product property for the invariant hyperbolic measures. \\

In this article, we focus on the holomorphic endomorphisms $f$ of
$\Pj^k$ of degree $d \geq 2$. These mappings define non invertible
ramified coverings of topological degree $d^k$. We refer to the
article of Dinh-Sibony \cite{DS} for a survey of their dynamical
properties. The question of the Hausdorff dimension for the equilibrium
measure was raised by Fornaess-Sibony \cite{FS2} (see subsection \ref{eqq}). \\

When $k=1$, $f$ defines a rational map on $\Pj^1$, and Ma\~n\'e
\cite{M} proved the formula $\delta = h_\nu /
\lambda$ $\pp$ for any ergodic measure satisfying $h_\nu > 0$. Here $\lambda$ denotes the single exponent of
$\nu$, it has multiplicity $2$ for the underlying real system. The
proof heavily relies on the Koebe distortion theorem. The present article deals with the higher dimensional
case, which is not conformal. We obtain the following result:

\begin{theorem}{A:} Let $f$ be a holomorphic endomorphism of $\Pj^k$
  of degree $d \geq 2$ and $\nu$ be an ergodic $f$-invariant
  measure with positive Lyapunov exponents $\lambda_k \leq \ldots \leq \lambda_1$. Then we have: 
\begin{equation*} \label{dere}
 \forall x \in \Pj^k \  \nu \textrm{-a.e.}  \ , \ \underline \delta(x)  \,  \geq \,  {\log d^{k-1} \over \lambda_1} + { h_\nu -  \log d^{k-1} \over \lambda_k } .\end{equation*} 
\end{theorem}

The proof is outlined in section \ref{pmk}, the method consists in
studying the distribution of the $\nu$-generic inverse
branches of $f^n$ in $\Pj^k$. Our main tools are a volume growth estimate for
holomorphic polydiscs in $\Pj^k$ and a normalization theorem for the
$\nu$-generic inverse branches of $f^n$. That result
provides, in some sense, a substitute for the one-dimensional Koebe distortion theorem.

\subsection{Application to the equilibrium measure $\mu$ of $f$}\label{eqq}

The equilibrium measure is defined as the limit (in the sense of distributions) of the smooth $(k,k)$ form
$d^{-kn} {f^n}^* \om^k$, where $\om^k$ is the standard volume form
on $\Pj^k$. Fornaess-Sibony \cite{FS1} proved that $\mu$ is mixing and
that $\log \Jac f \in L^1(\mu)$. Briend-Duval established that the
exponents of $\mu$ are bounded below by $\log \sqrt d$ \cite{BD1} and
that $\mu$ is the unique measure of maximal entropy ($h _\mu = \log d^k$) \cite{BD2}. Concerning the Hausdorff dimension of $\mu$, Ma\~n\'e's
formula asserts that $\dim_\HH \mu = \log d /
\lambda$ when $k=1$. Binder-DeMarco  \cite{BDeM} conjectured for $k \geq 2$:
\begin{conj}{:} For every system  $(\Pj^k,f,\mu)$, $\dim_\HH \mu  = { \log d \over \lambda_1 } + \cdots + { \log d \over \lambda_k }$. 
\end{conj}
We note that this formula is consistent with (\ref{ly}) if we set $h_i = \log d^i$ for the conditional entropies of $\mu$. Binder-DeMarco \cite{BDeM} proved  that $\dim_\HH \mu \leq 2k - 2
\, (\Sigma_{i=1}^k \lambda_i  -  k \log \sqrt d )/
\lambda_1$ in a polynomial setting by using volume estimates. Dinh-Dupont \cite{DD} extended that estimate to
meromorphic endomorphisms of $\Pj^k$. \\

From theorem A we deduce  the following bound. It proves half of the conjectured formula when $k=2$.
\begin{corollary}{A:} Let $f$ be an endomorphism of $\Pj^k$ of degree
  $d \geq 2$ and $\mu$ be its equilibrium measure. If $\lambda_k \leq  \cdots \leq \lambda_1$ denote the Lyapunov exponents of $\mu$, then 
\begin{equation*} \label{dere}
     \dim_\HH \mu   \,  \geq \,  {\log d^{k-1} \over \lambda_1} + { \log d \over \lambda_k } .
 \end{equation*} 
In particular, $\dim_\HH \mu    \geq   {\log d \over \lambda_1} + { \log d \over \lambda_2}$ for every system $(\Pj^2,f,\mu)$.
\end{corollary}
Now we can establish the conjecture for a class of non conformal systems by combining corollary A with the upper bound stated above:
\begin{corollary}{B:} Let $f$ be an endomorphism of $\Pj^k$ of degree
  $d \geq 2$ and $\mu$ be its equilibrium measure. If $\lambda_k = \log \sqrt d$ and $\lambda_{k-1} = \ldots = \lambda_1$, then
\[ \dim_\HH \mu = {\log d^{k-1} \over \lambda_1} + {\log d \over \lambda_k}. \] 
\end{corollary}

\subsection{Application to measures with large entropy}\label{evv}

Let $f$ be an endomorphism of $\Pj^k$ of degree $d \geq 2$ and $\nu$
  be an $f$-invariant ergodic measure. De Th\'elin \cite{dT} proved
  that if $\log \Jac f \in L^1(\nu)$ and  $h_\nu >
  \log d^{k-1}$, then the Lyapunov exponents of $\nu$ satisfy ${1 \over
  2}(h_\nu -  \log d^{k-1}) \leq \lambda_k \leq  \cdots \leq
  \lambda_1$. In \cite{Du} we recently constructed  ergodic measures satisfying $h_\nu >   \log d^{k-1}$ and showed that the preceding estimate holds without assuming the integrability of $\log \Jac f$. By theorem A, we deduce the following bounds for the largest Lyapunov exponent of $\nu$. 
    
\begin{corollary}{C:} Let $f$ be an endomorphism of $\Pj^k$ of degree $d \geq 2$ and $\nu$
  be an $f$-invariant ergodic measure.
 \begin{enumerate}
\item If  $\log d
  ^{k-1} < h_\nu$, then $\lambda_1 \geq  (1-1/k) \log \sqrt d$. 
\item If $\log d^{k-1} < h_\nu <  (1+
1/k) \log d^{k-1}$, then $\lambda_1 \geq {1 \over 2}(h_\nu -  \log
d^{k-1}) + \varphi(h_\nu)$, where $\varphi(h_\nu) > 0$. 
\end{enumerate}
\end{corollary}

The first point follows from the observation $\underline \delta \leq 2k$. For the second point, the function $\varphi$ is defined as $\varphi(h_\nu) := {1 \over 2} [ (1+1/k) \log
  d^{k-1} - h_\nu ]$. Let us observe that the latter is false for
 the equilibrium measure $\mu$, its Lyapunov exponents are indeed $\lambda_1 = \ldots = \lambda_k =  {1
  \over 2}(h_\mu -  \log d^{k-1}) = \log \sqrt d$ when $f$ is a Latt\`es example  \cite{BeDu}.

\subsection{Organization of the article}\label{orga}

The proof of theorem A relies on theorem B, which is stated in section 2: that result describes the distribution of the $\nu$-generic inverse
branches in $\Pj^k$. Section 3 deals with
notations and the normalization theorem for the inverse
branches. The proof of theorem A is detailed in section 4. We show theorem B in sections
5 and 6.  In an appendix we establish the growth lemma. \\

\emph{Acknowledgements :}  I thank the referee, whose advice and careful reading enabled me to improve the exposition of the article. Part of this work was written while visiting IMPA in Rio de Janeiro. I thank J.V. Pereira, M. Viana and the Institut for their kind hospitality.

\section{Statement of theorem B and outline of its proof}\label{pmk}
 
Let us fix $f$ a holomorphic endomorphism of $\Pj^k$ with degree $d \geq 2$ and $\nu$ an ergodic $f$-invariant measure with
positive exponents $\lambda_k \leq \ldots \leq \lambda_1$. The fractional time $q_n$ is defined as the entire part of $n \lambda_k
/ \lambda_1$. We denote by $f^{-n}_{y_n}$ the inverse branch of $f^n$
mapping $y_n := f^n(y)$
to $y$. We set $\EE_{\rho}$ as an arbitrary maximal $\rho$-separated subset in
$\Pj^k$. We define $\EE_{\rho} (q)$ as the finite set of $p \in \EE_{\rho}$ satisfying $q \in  B_p(\rho)$ and denote $B_x^{\Omega_\epsi}(r) := B_x(r) \cap \Omega_\epsi$. 

 \begin{theorem}{B :}  Let $f$ be an endomorphism of $\Pj^k$
  of degree $d \geq 2$ and $\nu$ be an ergodic $f$-invariant
  measure with positive Lyapunov exponents $\lambda_k \leq \ldots \leq
  \lambda_1$. For every $\epsi >0$, there exist $\Omega_\epsi
  \subset \Pj^k$ and $r_0 = r_0(\epsi) > 0$ satisfying: 
  \begin{enumerate}
\item $\nu(\Omega_\epsi) > 1 - \epsi$.
\item for every $x \in  \Omega_\epsi$ and $n$ large enough, the
  collection of inverse branches 
 \[ \PP_n(x) := \left \{  \, f^{-n}_{y_n} B_p (s_n) \ , \ y \in  B_x^{\Omega_\epsi} (s_n e^{-n
  \lambda_k + 3 n\epsi} )  \ , \ p \in  \EE_{s_n} (y_n) \,  \right \}    \]
 is well defined for $s_n := r_0 e^{-8n \epsi}$ and satisfies $\cdb \PP_n(x) \leq  d^{(k-1)(n - q_n)} e^{20kn\epsi}$.
\end{enumerate}
\end{theorem}

Theorem B is used in the proof of theorem A (see section \ref{uvb}). We sketch
below the proof of theorem B. It relies on propositions A and B. For simplicity, we shall work up to $e^{\pm n \epsi}$ error terms (for instance we replace $e^{-n
  \lambda_k + 3 n\epsi}$ by $e^{-n
  \lambda_k}$ and $s_n$ by $1$). \\

We define a \emph{polydisc} as any holomorphic map $\eta : \D^{k-1}(r) \to
\Pj^k$. Let $\om$ denote the Fubini-Study
$(1,1)$-form on $\Pj^k$ and define $\Vol \eta :=
\int_{\D^{k-1}(r)} \eta^* \om^{k-1}$: this is the \emph{volume} of $\eta$ counted with
multiplicity. Let $\{ B_j \, , \, j \in J
\}$ be a finite covering of $\Pj^k$ which consists of open sets bounded in the affine charts. We say that $\eta$ is \emph{bounded} if
its image is contained in some $B_j$. We shall need the
\begin{fact}{:}\label{aire}
If $\eta : \D^{k-1} (2) \to \Pj^k$ is bounded, then $\Vol  f^{m} \circ
\eta_{\vert \D^{k-1}}  \leq d^{(k-1) m}$ for every $m \geq 1$.
\end{fact}
That geometric result does not depend on the measure $\nu$: the proof  relies on the existence of a
Green current for every endomorphism of $\Pj^k$ (see the appendix). That lemma allows us
to establish the next proposition: let us fix $x \in
\Omega_\epsi$ and denote by $\LL_n$ the set of polydiscs
 $L_n :  \D^{k-1} \to B_x (e^{-n \lambda_k})$. 
 
\begin{proposition}{A :} For every $L_n \in \LL_n$, we have $\Vol f^n \circ L_n  \leq d^{(k-1)(n-q_n)}$.
\end{proposition}

This estimate follows from the growth lemma taking $m = n -
q_n$ and $\eta = f^{q_n} \circ L_n$. Indeed, the polydisc $f^{q_n} \circ L_n$
is bounded since $f^{q_n} ( B_x
(e^{-n \lambda_k}) ) \subset B_{x_{q_n}} (e^{-n \lambda_k} \cdot
e^{q_n \lambda_1}) \simeq B_{x_{q_n}}(1)$: that  comes from the fact that $\lambda_1$ is the largest exponent and $q_n \lambda_1  \simeq n \lambda_k$.\\

Our second tool is a normalization
theorem for the inverse branches of $f^n$ established by Berteloot-Dupont-Molino \cite{BDM}. That theorem basically asserts that every inverse branch $P_n \in \PP_n(x)$ looks like a parallelepiped with characteristic
dimensions $e^{-n \lambda_1} \leq \ldots \leq e^{-n \lambda_k}$, it plays
the role of a distortion theorem. The normalization theorem allows us to prove: 
\begin{proposition}{B :} There exists a finite subset $\FF_n \subset
  \LL_n$ of cardinality less than $k e^{20n\epsi}$ such that for every
  $P_n \in \PP_n(x)$, there is $L_n \in \FF_n$ satisfying $\Vol
  f^n \circ  {L_n}_{\vert L_n^{-1}(P_n)}  \geq 1$.
\end{proposition}
We actually show that $\Vol f^n \circ  {L_n}_{\vert
  L_n^{-1}(P_n)} \geq 1$ for (almost) every polydisc $L_n \in \LL_n$ transverse to the
$e^{-n\lambda_k}$-direction of $P_n$. The family $\FF_n$ then practically consists of hyperplanes parallel to the coordinates.\\

Finally, the upper bound $\cdb \PP_n(x) \leq  d^{(k-1)(n - q_n)}$ follows using the fact that the inverse branches are pairwise disjoint (see subsection \ref{cedec}), that completes the proof of theorem B.\\

Let us notice that the estimates of theorems A and B can be sharpened when
$\lambda_k$ has multiplicity $p$. The same method indeed yields $ \cdb
\PP_n(x) \leq d^{(k-p)(n - q_n)}$ by considering the family of
polydiscs $L_n :  \D^{k-p} \to B_x (e^{-n \lambda_k})$. In particular
that implies the lower bound $\dim_\HH \nu \geq {\log d^{k-p} \over \lambda_1 } + {h_\nu - \log d^{k-p} \over \lambda_k}$.

\section{Generalities}\label{zpp}

\subsection{The dynamical systems $(\Pj^k,f,\nu)$}\label{ds}

Let $f$ be a holomorphic endomorphism of $\Pj^k$ of degree $d \geq
2$. It is defined in homogeneous coordinates as $[ P_0 :  \ldots : P_k ]$ where the $P_i$'s are homogeneous
polynomials of degree $d$ without common zero (except the origin). The topological degree
of $f$ is $d^k$. Let $\CC$ be the
critical set of $f$, this is an hypersurface of degree
$(d-1)(k+1)$ counted with multiplicity. Let $\om$ be the Fubini-Study
$(1,1)$-form on $\Pj^k$ and $\dist$ the induced distance on $\Pj^k$. We denote by
$\Jac f$ the function on $\Pj^k$ satisfying $f^* \om^k = \Jac f
\cdot \om^k$. This is a bounded non-negative $C^\infty$ function which vanishes on the
critical set $\CC$.\\

Let $\nu$ be an $f$-invariant ergodic measure, $h_\nu$ its entropy and $\lambda_k \leq \ldots \leq
  \lambda_1$ its Lyapunov exponents. We assume that those exponents are
  positive. In particular, the
  classical formula $\int_{\Pj^k} \log \Jac f \, d\nu = 2 (\lambda_1 +
  \ldots +\lambda_k )$ yields:

\begin{lem}\label{integ}
If the exponents of $\nu$ are positive, then $\log \Jac f \in
L^1(\nu)$ and $\nu(\CC) = 0$. 
\end{lem}

We shall assume that $\lambda_k < \ldots < \lambda_1$. In particular that enables us to simplify the statements concerning the
normal forms (see the next subsections). Our method easily extends
when multiplicities occur. \\

We endow $\C^k$ with $\abs{z} = \max_{1 \leq i \leq k}  \abs{z_i}$. For any polynomial
mapping $Q : \C^k \to
\C^l$, we set $\norm{Q}$ as the maximum
of the modulus of its coefficients. We also denote by $(c_i)_{1 \leq i \leq
  k}$ the canonical basis of $\C^k$ and by $(\pi_i)_{1 \leq i \leq k}$ the projections to the axis.

\subsection{Normal forms associated with the Lyapunov exponents}\label{lkj}

For every $\al =
(\al_1,\ldots,\al_k) \in \N^k$, we set $\abs{\al} := \al_1 + \ldots + \al_k$ and $Q_\al := z_1^{\al_1} \ldots
z_k^{\al_k}$.
Given $1 \leq i \leq k-1$, the set of $i$-resonant degrees is defined by: 
\[ \RRR_i := \left \{ \al \in \N^k \ , \  \abs{\al} \geq 2 \ , \ \al_1 =
\ldots = \al_i = 0 \  \textrm{ and }  \  \lambda_i =  \al_{i+1} \lambda_{i+1}
  +  \ldots + \al_k  \lambda_k  \right \}. \] 
We set $I := \{ 1 \leq i \leq k-1 \, , \, 2 \lambda_k \leq \lambda_i  \}$. Observe that $\RRR_i$ is empty if $i \notin I$. Note also that $\abs{\al} \leq \theta := \lambda_1 / \lambda_k$ for every $\al
  \in \RRR_i$, hence $\RRR := \cup_{i=1}^{k-1} \RRR_i$ has finite cardinality. We denote $\Delta := \cdb \RRR$.\\

We say that a polynomial map $N : \C^k \to \C^k$ is \emph{normal} if
 $N = ( N_1 , \ldots , N_{k-1},0)$ where $N_i =
\sum_{\al \in \RRR_i} c_i^\al Q_\al$ for some $c_i^\al \in \C$. A map $R : \C^k \to \C^k$ is \emph{resonant} if $R = A+N$, where  $A = (a_1,\ldots,a_k)$ is a linear diagonal map satisfying 
 $ e^{-\lambda_i - \epsi} \leq \abs{a_i} \leq e^{-\lambda_i + \epsi}$ and $N$ is a normal
map. \\

Every resonant map $R = A+N$ is invertible, and $R^{-1} = A^{-1} + N'$ for some normal map $N'$. Moreover, if $R_i = A_i + N_i$ ($i =
1,2$) are resonant maps, we have $R_1 \circ R_2 = A_1 \circ A_2 + N''$
for some normal map $N''$. These are classical stability properties
(see e.g. \cite{GK}, section 1.1 and \cite{BDM}, section 5).

\subsection{Natural extension and normalization theorem}\label{kop}

Let $\orb := \{ \hat x := (x_n)_{n \in \Z} \, , \,  x_{n+1} = f (x_n) \}$ be the set of orbits, $\hat \pi : \orb \to \Pj^k$ the projection $\hat x \mapsto
x_0$, and $s : \orb \to \orb$  the left shift. We also set $\tau :=
s^{-1}$. Note that $\hat \pi \circ s = f \circ \hat
\pi$ on $\orb$. For every $n \geq 0$,
we denote $\hat x_n := s^n(\hat x)$. We say that a function
$\phi_\epsi : \orb \to \R^+$ is
	\emph{$\epsi$-slow} (resp. \emph{$\epsi$-fast}) if $\phi_\epsi(\orb) \subset \, ]0,1]$ (resp. $[1,+\infty[$) and satisfies
		$\phi_\epsi (\hat
		x) e^{- \epsi}  \leq \phi_\epsi(s(\hat x)) \leq \phi_\epsi (\hat
		x) e^{\epsi}$ for every $\hat x \in \orb$. \\
		
We denote by $\hat \nu$ the $s$-invariant measure on $\orb$ satisfying $\hat \nu (\hat \pi ^{-1}(A)) = \nu(A)$ for
  every borel set $A \subset \Pj^k$ (see \cite{CFS}, section
  10.4). We shall work with the $s$-invariant set $X :=  \{ \hat x = (x_n)_{n \in
  \Z} \, , \, x_n \notin \CC \}$. It satisfies $\hat \nu(X) = 1$
  since $\nu(\CC) = 0$ (see lemma \ref{integ}). For every $\hat
x \in X$, we denote by $f^{-n}_{\hat x}$ the inverse branch of $f^n$
  sending $x_0$ to $x_{-n}$. Hence $f^{-n}_{\hat x_n}$ is the inverse
  branch of $f^n$ sending $x_n = f^n(x)$ to $x$.

\begin{defn}
$\RR = (R_{\hat x})_{\hat x \in
  X}$ is a \emph{resonant cocycle} if every $R_{\hat x}$ is a resonant
  map. 
\end{defn}
  
Given a resonant cocycle $\RR$, we set $R_{\hat x} := \left( a_{1}({\hat
  x}),\ldots,a_{k}({\hat x}) \right)  + \left( N_{1}(\hat x),
  \ldots, N_{k-1}(\hat x) , 0  \right)$, where  $ e^{-\lambda_i -
  \epsi} \leq \abs{a_i(\hat x)}
  \leq e^{-\lambda_i + \epsi}$. For every $n \geq 1$,
  we define $R^n_{\hat x} := R_{\tau^{n-1}(\hat x)}  \circ \ldots \circ R_{\hat
  x}$ and $R^{-n}_{\hat x} := (R^{n}_{\hat x})^{-1}$. Using the stability
  properties, we obtain:
\begin{equation}\label{ztop}
\forall n \in \Z \ , \ R^n_{\hat x}  = \left( a_{1,n}({\hat
  x}),\ldots,a_{k,n}({\hat x}) \right)  + \left( N_{1,n}(\hat x),
  \ldots, N_{k-1,n}(\hat x) , 0  \right) ,
\end{equation}
where $e^{-n \lambda_i - \abs{n}
  \epsi} \leq \abs{a_{i,n}(\hat x)} \leq e^{-n \lambda_i + \abs{n}
  \epsi}$ and $N_{i,n}(\hat x) := \sum_{\al \in \RRR_i}
  c_{i,n}^\al(\hat x) Q_\al$.
  
 \begin{defn} Let $M_\epsi$ be an $\epsi$-fast function on $X$. A resonant cocycle $\RR$ is \emph{$M_\epsi$-adapted} if $\norm{ N_{i,n}(\hat x) } = \max_{\al \in \RRR_i} \abs{c_{i,n}^\al(\hat x)} \leq
M_\epsi(\hat x) e^{-n\lambda_i + \abs{n}  \epsi}$ for every $n \in \Z$. 
\end{defn}

\begin{defn}\label{ssss} Let $r_\epsi$ and $\beta_\epsi$ be
  respectively an $\epsi$-slow and an
  $\epsi$-fast function on $X$. $\SS = (S_{\hat x})_{\hat x \in X}$ is a \emph{($r_\epsi,\beta_\epsi$)-coordinate} if for any $\hat x \in X$,  $S_{\hat x} : B_{x_0}(r_\epsi(\hat x)) \to
  \C^k$ is an injective holomorphic map satisfying $S_{\hat x}(x_0)=0$ and
\begin{equation*}\label{klm}
 \forall (p,p') \in B_{x_0}(r_\epsi(\hat x)) \ , \ \dist(p,p') \leq
  \abs{S_{\hat x}(p) - S_{\hat x}(p') } \leq  \beta_\epsi(\hat x) \,
  \dist(p,p') .
\end{equation*}
\end{defn}

The normalization theorem  is stated as follows \cite{BDM}.

\begin{thm}\label{nt}
For every $\epsi >0$, there exist a
$(r_\epsi,\beta_\epsi)$-coordinate $\SS$ and an $M_\epsi$-adapted resonant
cocycle $\RR$ such that the following diagram commutes for $\hat
\nu$-almost every $\hat x \in X$ and every $n \geq 1$:
\[ 
\xymatrix{
  B_{x_0}(r_\epsi(\hat x)) \ar[rr]^{f^{-n}_{\hat x}}  \ar[d]_{ S_{\hat x}}  &   &  f^{-n}_{\hat x}  \left( B_{x_0} (r_\epsi(\hat x)) \right) \ar[d]^{S_{\tau^n(\hat x)}}   \\
     \C^k   \ar[rr]^{R^n_{\hat x}}         &   &  \C^k   } \]
\end{thm}

Note that the existence of $r_\epsi$ requires the $\nu$-integrability of $\log \norm{ (d_x
  f)^{-1}}$ (see \cite{BDM}, lemma 4.1). Here this is a consequence of
  lemma \ref{integ}. 

\subsection{Some estimates}\label{ggm}

We denote $z := (\tilde z , z_k) \in \D^{k-1}
\times \D$ and $\tilde \pi (z) := \tilde z$. We recall that $\Delta =
\cdb \RRR$ and that $\abs{\al} \leq \theta = \lambda_1 / \lambda_k$ for every $\al \in
\RRR$.

\begin{lem}\label{bnormb} Let $\RR$ be an $M_\epsi$-adapted resonant cocycle
  and $M'_\epsi :=  \max \{ \Delta + 1 , \theta , \theta(\theta -1) \}
  M_\epsi$. Then for every $\hat x \in X$, $r \leq 1$ and $z  \in \D^k(r)$,
  we have:
\begin{enumerate} 
\item $\D^k \left(  M'_\epsi (\hat x)^{-1} e^{-n\lambda_1 - n \epsi}
  \cdot r \right) \subset R^n_{\hat
  x}\left( \D^k(r) \right) \subset \D^k \left( M'_\epsi(\hat x)
  e^{-n\lambda_k + n\epsi} \cdot r \right)$,
 
\item $ \norm { \tilde \pi \circ d_z  R^n_{\hat x} }  \leq M'_\epsi(\hat x) e^{-n\lambda_{k-1} + n\epsi}$,

\item $ e^{-n \lambda_k -  n\epsi} \leq  \left \vert \pi_k  \left (
  {\partial  R^n_{\hat x} \over \partial z_k}  (z) \right ) \right \vert $
  \, and \,  $\left \vert {\partial  R^n_{\hat x} \over \partial z_k} (z)  \right \vert \leq \max \{ M'_\epsi(\hat x) e^{-n\lambda_{k-1} + n\epsi}   \, , \,    e^{-n \lambda_k + n\epsi} \}$,

\item $\left \vert {\partial^2  R^n_{\hat x} \over \partial z_k^2} 
  (z) \right \vert \leq  M'_\epsi(\hat x) e^{-2n\lambda_k +  n \epsi}$.  
\end{enumerate}
\end{lem}

\proof Let $\hat x \in X$, $r \leq 1$ and $z  \in \D^k(r)$. Using the $M_\epsi$-adapted property and (\ref{ztop}), we get for every $1 \leq i \leq k$:
\[ \abs{ \pi_i (R^n_{\hat x} (z) )  } \leq \abs{
  a_{n,i}(\hat x) } \abs{z} +  \Delta \, \norm { N_{i,n}(\hat x)  } \abs{z}^{\theta} \leq
  (\Delta + 1) M_\epsi(\hat x) e^{-n\lambda_i + n\epsi} \abs{z} .\]
We deduce $\abs{ R^n_{\hat x} (z) } < M'_\epsi(\hat x) e^{-n\lambda_k +
  n\epsi} r$ for every $z \in \D^k(r)$. Similarly, for every $w
\in \D^k(r)$ and  $1 \leq i \leq k$, we have:
\[ \abs{ \pi_i (R^{-n}_{\hat x} (w) )  } \leq (\Delta + 1) M_\epsi(\hat x) e^{n\lambda_i +
  n\epsi} \abs{w} \leq  M'_\epsi(\hat x) e^{n\lambda_1 + n\epsi} \abs{w}. \]
Hence $\abs{ R^{-n}_{\hat x} (w)} < r$ for every $w \in \D^k \left(  M'_\epsi (\hat x)^{-1} e^{-n\lambda_1 - n \epsi}
  r \right)$. That proves the point 1. For the point 2, observe that for every $1 \leq
  i \leq k-1$ and $z \in \D^k(r)$:
\[ \norm{ \pi_i \circ d_z R^n_{\hat x}  } \leq \max \{ \, \abs{
  a_{n,i}(\hat x) } \, , \, \theta \norm{ N_{i,n}(\hat x)  } r^{\theta - 1} \,   \}
  \leq M'_\epsi(\hat x) e^{-n\lambda_{k-1} + n\epsi} . \]
The point 3 now follows from the point 2 and the observation (see
  (\ref{ztop})) : 
\[   \left \vert \pi_k  \left ( {\partial  R^n_{\hat x} \over \partial z_k} (z) \right ) \right \vert  =   \norm{\pi_k \circ d_z R^n_{\hat x}} = \abs{a_{k,n} (\hat
  x)} \simeq e^{-n\lambda_k \pm n\epsi}. \]
For the point 4, let us distinguish whether or not $I = \{ 2 \lambda_k \leq \lambda_i  \}$ is empty. If $I$ is empty, there are no resonant degree, hence $R^n_{\hat x}$ is a linear mapping and ${\partial^2
  R^n_{\hat x} \over \partial z_k^2} = 0$. If $I$ is not empty ($\theta = \lambda_1 / \lambda_k \geq 2$ in that case), we have for every $1 \leq i \leq \max I$:
\[  \left \vert   \pi_i \left(  {\partial^2  R^n_{\hat
  x} \over \partial z_k^2} (z) \right) \right \vert \leq \theta(\theta-1) \norm{ N_{i,n}(\hat x) } r^{\theta-2}
  \leq  M'_\epsi(\hat x)
  e^{-n\lambda_i + n\epsi} \leq  M'_\epsi(\hat x) e^{-2n\lambda_k +
  n\epsi}, \] 
and $\pi_i ( {\partial^2  R^n_{\hat
  x} \over \partial z_k^2} ) = 0$ for every $\max I +1 \leq  i \leq
  k$.  \finsec

\section{Proof of theorem A}\label{uvb}

In this section we establish theorem A assuming theorem B. Our aim is to prove:
\begin{equation} \label{dere}
  \forall x \in \Pj^k \  \nu\textrm{-a.e.} \ , \ \liminf_{r \to 0} \, { \log \nu (B_x(r)) \over  \log r } \, \geq \,
   {\log d^{k-1} \over \lambda_1} + { h_\nu - \log d^{k-1}  \over \lambda_k } .
\end{equation} 
Let $\epsi > 0$ and $\Omega_\epsi$, $r_0$ be given by theorem B. We have $\nu(\Omega_\epsi) > 1-\epsi$
, and for every $x \in \Omega_\epsi$ the cardinality of  
\[ \PP_n(x) =  \{ \,  f^{-n}_{\hat
  y_n}(B_p(s_n)) \, , \, y \in B_x^{\Omega_\epsi}(\rho_n)  \, , \,  p
\in \EE_{s_n}(y_n) \, \} \]
is less than $d^{(k-1)(n - q_n)} e^{20kn\epsi}$. Here we set $\rho_n := s_n e^{-n \lambda_k
  + 3n\epsi}$. We shall use
Brin-Katok's theorem. Let $B_n(x,\xi) := \{  z \in \Pj^k \, , \, \dist (f^q(x) , f^q(z))
<  \xi \, , \,  0 \leq q \leq n \}$ be the $n$-dynamical ball centered at $x$ with radius $\xi$.
\begin{theorem}{\cite{BK}} For $\nu$-a.e. $x \in \Pj^k$, we have 
\[ \sup_{\xi > 0} \,  \liminf_{n \to +\infty} \,  - {1 \over n}
\log \nu (B_n(x,\xi)) = h_\nu. \]
\end{theorem}
\noindent In particular, for $\nu$-a.e. $x \in \Pj^k$, there exist $\xi_\epsi(x) > 0$ and $m_\epsi(x) \geq 1$ such that:
\begin{equation*}\label{entbk}
\forall  \xi \leq  \xi_\epsi(x) \, , \, \forall  n \geq m_\epsi(x) \, , \, \nu(B_n(x,\xi)) \leq e^{-n(h_\nu-\epsi)} .
\end{equation*}
We may decrease $r_0$ and choose $m_0 \geq 1$ large enough so that $\Gamma_\epsi := \{ \xi_\epsi
\geq r_0 \, , \, m_\epsi  \leq m_0  \} $ satisfies
$\nu(\Gamma_\epsi) > 1-\epsi$. We have:
\begin{equation}\label{entbkkk}
 \forall x \in \Gamma_\epsi \, , \, \forall n \geq m_0 \, , \, \nu(B_n(x,r_0)) \leq e^{-n(h_\nu-\epsi)} .
\end{equation}
 We let $\Lambda_\epsi := \Gamma_\epsi \cap
\Omega_\epsi$ (it satisfies $\nu(\Lambda_\epsi) > 1-2\epsi$) and define:
 \[ \QQ_n(x) :=  \{ \,  f^{-n}_{\hat
  y_n}(B_p(s_n)) \, , \, y \in B_x^{\Lambda_\epsi}(\rho_n)  \, , \,  p
 \in \EE_{s_n}(y_n) \, \} \ \subset \ \PP_n(x). \]

\begin{lem}\label{pmlk}
For every $Q \in \QQ_n(x)$ we have $\nu ( Q ) \leq  e^{-n(h_\nu -\epsi)}$.
\end{lem}

The proof needs the definition of $\Omega_\epsi$ and is postponed to subsection \ref{defam}. Let $\Lambda_\epsi' \subset \Lambda_\epsi$ be the subset of points satisfying $\nu (B_x^{\Lambda_\epsi} (r) ) / \nu (B_x(r)) \to 1$ when
$r \to 0$. The Borel density lemma asserts that $\nu(\Lambda_\epsi') = \nu(\Lambda_\epsi)$.

\begin{lem}\label{where2} For every $x \in \Lambda'_\epsi$, there exists $p(x) \geq 1$ such that: 
\[ \forall n \geq p(x) \ , \ \nu \left( B_x(\rho_n) \right) \leq 2 \, \cdb \PP_n (x) \cdot e^{-n(h_\nu-\epsi)}. \]
\end{lem}

\proof  Let $x \in \Lambda'_\epsi$ and $p(x) \geq 1$ so that $\nu
(B_x(\rho_n)) \leq 2 \nu (B_x^{\Lambda_\epsi} (\rho_n))$ for  $n
\geq p(x)$. The fact that $\QQ_n(x)$ is a covering of $B_x^{\Lambda_\epsi}(\rho_n)$ combined with lemma \ref{pmlk} implies: 
\[  \nu (B_x^{\Lambda_\epsi}(\rho_n))  \leq  \sum_{Q \in \QQ_n(x)} \nu( Q ) \leq  \, \cdb
\QQ_n (x) \cdot e^{-n(h_\nu-\epsi)}.  \] 
We conclude using $\cdb \QQ_n(x) \leq \cdb \PP_n(x)$. \fin

\begin{lem}\label{cucu}  For every  $x \in \Lambda_\epsi'$, we have:
\begin{equation*} \label{dff}
 \liminf_{r \to 0} \, { \log \nu (B_x(r)) \over  \log r } \, \geq \, \left(
   { \log d^{k-1} \over \lambda_1 }  + {h_\nu - \log d^{k-1} \over
   \lambda_k}  - {21 k \epsi \over \lambda_k } \right) {\lambda_k \over \lambda_k + 5\epsi}.
\end{equation*} 
\end{lem}

\proof  Lemma \ref{where2} yields for $n \geq p(x)$: 
\begin{equation*}
 \log \nu (B_x(\rho_n))   \leq   \log \cdb \PP_n(x) - n (h_\nu-\epsi)  + \log 2 .
\end{equation*}
We use theorem B to obtain for $n \geq p(x)$:
\begin{equation*}
 \log \nu (B_x(\rho_n))   \leq  (n-q_n) \log d^{k-1} - n (h_\nu - \epsi) +  20 k n \epsi + \log 2 .
\end{equation*} 
Using $\rho_n = r_0 e^{-n \lambda_k - 5n\epsi}$ and $q_n \geq n
\lambda_k / \lambda_1 - 1$, we obtain for $n$ large:
\[  { \log \nu (B_x(\rho_n)) \over \log \rho_n } \geq { 
 n \lambda_k / \lambda_1  \cdot \log d^{k-1}   + n ( h_\nu - \log d^{k-1} ) - 21 k n \epsi -  \log 2 \over n\lambda_k +  5n\epsi - \log r_0  }. \]
The aimed estimate follows  letting $n \to \infty$. \fin
 
Finally, lemma \ref{cucu} yields (\ref{dere}) as follows. Let $\Lambda' :=
\cap_{p \geq 1} \cup_{q \geq p} \Lambda_{1/q}'$. We have
$\nu(\Lambda') = 1$ since $\nu(\Lambda_{1/q}') > 1 - 2/q$ for every $q
\geq 1$. Now for every $x \in \Lambda'$ there exists a subsequence $(q_j(x))_{j \geq 1}$ such that
$x \in \Lambda_{1/q_j(x)}'$. We deduce (\ref{dere}) from lemma \ref{cucu} setting $\epsi =
1/q_j(x)$ and letting $j \to \infty$. That completes the proof of theorem A.

\section{Proof of theorem B}\label{defom}

\subsection{Definition of $\Omega_\epsi$ and $r_0$}\label{defam}

Let $\epsi > 0$ and $r_\epsi$, $\beta_\epsi$, $M'_\epsi$ be the $\epsi$-slow and $\epsi$-fast functions provided by theorem \ref{nt} and lemma \ref{bnormb}. Let us choose $r_0 \leq  1$ small 
and $\beta_0 , M_0' \geq 1$ large such that the set
\[ \widehat \Omega_\epsi :=
\left \{ \, \hat x \in X \, , \, r_\epsi(\hat x) \geq r_0 \,
, \, \beta_\epsi(\hat x) \leq \beta_0  \, , \, M'_\epsi(\hat x) \leq M_0' \,  \right\} \]
satisfies $\hat \nu
(\widehat \Omega_\epsi) > 1 - \epsi$. We define $\Omega_\epsi := \hat \pi(\widehat \Omega_\epsi)$. Observe that $\nu(\Omega_\epsi) = \hat \nu (\hat \pi^{-1}(\Omega_\epsi)) \geq
\hat \nu (\widehat \Omega_\epsi) > 1 - \epsi$.  We fix once and for all a section of
the restriction $\hat \pi :  \widehat  \Omega_\epsi \to \Omega_\epsi$. That is to say that we associate to every $x \in \Omega_\epsi$ an element
of $\widehat \Omega_\epsi \cap  \hat \pi^{-1} \{ x \}$, that
we denote $\hat x$.\\

We set $r_n := r_0 e^{-n\epsi}$, $\beta_n := \beta_0
e^{n\epsi}$, $M_n' := M_0'
e^{n\epsi}$. We shall also need:
\[  s_n :=  r_0 e^{-8 n\epsi} \ \ , \ \ \rho_n := r_0 e^{-n \lambda_k
  - 5n\epsi} \ \ , \ \  \tau_n :=  \beta_0 (1+ 2\beta_0 M_0' ) \rho_n. \]
In the sequel, the estimates and inclusions will be written for $n$
  large only depending on $\epsi , r_0, \beta_0 , M_0'$ and $(\lambda_i)_{1 \leq i
  \leq k}$.

\begin{lem}\label{nouv} For every $x \in \Omega_\epsi$, the maps $f^{-n}_{\hat x_n}$, $S_{\hat x_n}$ and  $R^n_{\hat x_n}$ satisfy:
\begin{enumerate} 
\item $f^{-n}_{\hat x_n}$ and $S_{\hat x_n}$ are well defined on
$B_{x_n}(r_n)$. 

\item  $\dist(p,p') \leq \abs{S_{\hat x_n}(p) - S_{\hat x_n}(p') } \leq  \beta_n  \, \dist(p,p')$ for every $(p,p') \in B_{x_n}(r_n)$. 

\item $\D_{ S_{\hat x_n} (p) } ^k(r) \subset S_{\hat x_n}(B_p(r)) \subset
\D_{S_{\hat x_n} (p)} ^k(\beta_n r)$ for every $B_p(r) \subset B_{x_n}(r_n)$.

\item $\D^k \left(  {M_n'}^{-1} e^{-q\lambda_1 - q \epsi}
  \cdot r \right) \subset R^q_{\hat
  x_n}\left( \D^k(r) \right) \subset \D^k \left( M_n' e^{-q\lambda_k +
  q\epsi} \cdot r \right)$ for every $r \leq 1$ and $0 \leq q \leq n$.

\end{enumerate} 
\end{lem}

\proof The fact that $\hat x \in  \widehat \Omega_\epsi$ and the $\epsi$-slow, $\epsi$-fast properties of $r_\epsi , \beta_\epsi$ yield $r_\epsi(\hat x_n) \geq r_\epsi(\hat x)
e^{-n\epsi} \geq r_n$ and $\beta_\epsi(\hat x_n) \leq
\beta_\epsi(\hat x) e^{n\epsi} \leq \beta_n$. All the items then
follow from  theorem \ref{nt}, definition \ref{ssss} and lemma \ref{bnormb}(1). \fin 

Now we can give the \proo\textsc{of lemma \ref{pmlk}:}  Let $y \in \Lambda_\epsi$ and $p \in \EE_{s_n}(y_n)$ such that
  $Q = f^{-n}_{\hat
  y_n}(B_p(s_n))$. Since $B_{p}(s_n) \subset B_{y_n}(2s_n)$, it
  suffices to prove that $f^{-n}_{\hat y_n}(B_{y_n}(2s_n))
\subset B_n (y,r_0)$ (see (\ref{entbkkk})). We verify for that purpose that $\dist(f^{-q}_{\hat
  y_n}(z) , f^{-q}_{\hat y_n}(y_n)) \leq r_0$  for every $z \in
  B_{y_n}(2s_n)$ and $0 \leq q \leq n$. Using the identity $f^{-q}_{\hat
  y_n} = S_{\hat y_{n-q}}^{-1} \circ  R^q_{\hat y_n} \circ S_{\hat
  y_n}$ and lemma \ref{nouv}(3,4), we get:
\[ \forall z \in B_{y_n}(2s_n) \ , \  \dist (f^{-q}_{\hat y_n}(z) ,
  f^{-q}_{\hat y_n}(y_n)) \leq 2 s_n \beta_n M_n' e^{-q \lambda_k + q
    \epsi}  \leq 2 r_0 M_0' \beta_0 e^{-5 n \epsi} \leq r_0.    \]
That completes the proof of lemma \ref{pmlk}. \fin

Let us deal with the biholomorphism $\psi_{x,y} := S_{\hat x} \circ S_{\hat y}^{-1}$ when $y$ is close to $x \in \Omega_\epsi$.

\begin{lem}\label{controun}
There exist $R \leq 1$ and $\ga > 0$ such that for every
$x \in \Omega_\epsi$ and $y \in B_x^{\Om_\epsi} (r_0/2)$:  
\begin{enumerate}
\item $\psi_{x,y}  : \D^k(R) \to \D^k(\beta_0)$ is well defined.

\item   ${1 \over \beta_0} \abs{z-z'} \leq \abs{\psi_{x,y} (z) - \psi_{x,y}
  (z')  } \leq  \beta_0  \abs{z-z'}$ for every $(z,z') \in \D^k(R)$.

\item $\norm{ d_z \psi_{x,y}  - d_{z'} \psi_{x,y} } \leq
 \ga \abs{ z-z'  }$ for every $(z,z') \in \D^k(R)$.
\end{enumerate}
\end{lem}
The point 2 actually implies $\abs{ d_0 \psi_{x,y}(c_k)  }
\geq 1/\beta_0$. We therefore have $B_x^{\Om_\epsi} (r_0/2) = \cup_{i=1}^k
W_x^i$, where $W_x^i := \{ \, \abs{ \, \pi_i \, ( d_0 \psi_{x,y}(c_k) ) \, } \geq 1 / \beta_0
      \}$. We fix for every $x \in \Omega_\epsi$ a partition $B_x^{\Om_\epsi} (r_0/2) = \cup_{i=1}^k Y_x^i$, where $Y_x^i
  \subset W_x^i$. We complete lemma \ref{controun} as follows: 
\begin{lem}\label{ggp}
$\abs{\pi_i \left(
  d_z \psi_{x,y} (c_k) \right)} \geq  1/ (2\beta_0)$ for every $y \in Y^i_x$ and $z \in \D^k(R)$. 
\end{lem}

\noindent \textsc{Proof of lemmas \ref{controun} and \ref{ggp}:} Let $R' = r_0 / 2$, $\ga = \beta_0 /
R'^2$ and $R = 1/ (2\beta_0 \ga) < R'$. We prove 1, 2 on $\D^k(R')$ and 3, lemma \ref{ggp} on
$\D^k(R)$. Lemma \ref{nouv}(3) yields for every $w \in \{ x , y \}$ (take $p = w$ and $n = 0$ in that lemma):
\begin{equation}\label{cb}
 \forall r \leq r_0  \ , \  \D^k (r) \subset S_{\hat w} (B_w (r)) \subset
\D^k (\beta_0 r).
\end{equation}
Let $z  \in \D^k(R')$. The left inclusion in (\ref{cb}) with $w=y$ and $r=R'$
yields $S_{\hat y}^{-1} (z) \in B_y(R')$.  Since $B_y(R') \subset B_x(r_0)$, the right inclusion in (\ref{cb}) with $w =
x$ gives $\psi_{x,y}(z)  = S_{\hat x} \circ S_{\hat y}^{-1} (z) \in
\D^k(r_0\beta_0) \subset \D^k(\beta_0)$. That proves the point 1. The
point 2 then comes from lemma \ref{nouv}(2) and the point 3 from
Cauchy's estimates: we indeed have $\norm{\psi_{x,y}}_{ C^2, \D^k(R')} \leq \beta_0 / {R'}^2 =
\ga$ from point 1. Now let us deal with lemma \ref{ggp}. For every $z \in
\D^k(R)$, the point 3
implies  $\norm{ d_z \psi_{x,y} - d_0
  \psi_{x,y} } \leq \ga R = 1 / (2\beta_0)$. The desired estimate then
follows from $\abs{ \pi_i ( d_0 \psi_{x,y} (c_k) )}
\geq 1/  \beta_0$. \finsec

\subsection{The upper bound on $\cdb \PP_n(x)$}\label{cedec}

Let $x \in \Omega_\epsi$. Recall that $s_n = r_0 e^{-8n\epsi}$, $\rho_n = s_n e^{-n
  \lambda_k +3n\epsi}$ and 
\begin{equation*} 
  \PP_n(x) = \left \{ \, f^{-n}_{\hat y_n}( B_p(s_n) ) \, , \,  y \in B_x^{\Om_\epsi}(\rho_n) \, , \,  p \in \EE_{s_n}(y_n)  \right \},
\end{equation*} 
where $\EE_{s_n}$ is a fixed $s_n$-separated
  set in $\Pj^k$. We want to prove 
\begin{equation}\label{tha}
\cdb \PP_n(x) \leq   d^{(k-1)(n-q_n)} \cdot e^{ 20 k n\epsi},
\end{equation} 
where $q_n$ denotes the entire part of $n \lambda_k/\lambda_1$. We first verify
that $\PP_n(x)$ is well defined, it therefore induces a covering of
$B_x^{\Om_\epsi}(\rho_n)$:

\begin{lem}\label{gds} For every $y \in
\Om_\epsi$ and $p \in \EE_{s_n}(y_n)$:
\begin{enumerate}
\item $f^{-n}_{\hat y_n}$  and $S_{\hat y_n}$ are well defined on $B_p(s_n)$,
\item $S_{\hat y_n} (B_p(s_n)) \subset \D^k(2 s_n \beta_n)$.
\end{enumerate}
\end{lem}

\proof Observe that $B_p(s_n) \subset B_{y_n}(2 s_n) \subset
B_{y_n}(r_n)$ by definition of $\EE_{s_n}(y_n)$. The items then
follows from lemma
\ref{nouv}(1,3). \fin

Now we localize the collection $\PP_n(x)$. We recall that $\tau_n
= \beta_0(1+2 \beta_0 M_0') \rho_n$.
\begin{lem}\label{where}  For every $x \in
\Om_\epsi$ and $P_n \in \PP_n(x)$, we have $P_n \subset S_{\hat x}^{-1} ( \D^k(\tau_n) )$. 
\end{lem}

\proof Let $P_n \in \PP_n(x)$: there exist $y \in
B_x^{\Om_\epsi}(\rho_n)$ and $p \in \EE_{s_n}(y)$ satisfying $P_n =
f^{-n}_{\hat y_n}(  B_p(s_n) )$. Our aim is to prove that
$S_{\hat x}(P_n) \subset \D^k(\tau_n)$. We shall use $S_{\hat
  x}(P_n) = \psi_{x,y}  \circ  R^n_{\hat y_n} \circ S_{\hat y_n}(P)$, where $P :=  B_p(s_n)$. This comes from $f^{-n}_{\hat y_n} =
S_{\hat y}^{-1} \circ  R^n_{\hat y_n} \circ S_{\hat y_n}$ (see theorem
\ref{nt}) and $\psi_{x,y} = S_{\hat x} \circ S_{\hat y}^{-1}$. Lemmas
\ref{gds}(2) and \ref{nouv}(4) yield successively $S_{\hat y_n}(P) \subset
\D^k(2 s_n \beta_n)$ and $R^n_{\hat y_n} \circ S_{\hat y_n}(P) \subset
\D^k(2 s_n \beta_n M_n' e^{-n\lambda_k + n\epsi})$, which is included in $\D^k(R)$. Then lemma \ref{controun}(1,2) implies:
\begin{equation}\label{gty}
 \psi_{x,y} \circ  R^n_{\hat y_n} \circ S_{\hat y_n}(P)  \subset
 \psi_{x,y}(0) + \D^k(2 s_n  \beta_n M_n' e^{-n\lambda_k + n\epsi}
 \beta_0). 
\end{equation} 
But $\psi_{x,y} (0) = S_{\hat x}(y) \in \D^k( \rho_n \beta_0)$ from $y \in  B_x(\rho_n)$ and lemma \ref{nouv}(3). The
  right hand side of (\ref{gty}) is therefore included in $\D^k(\rho_n \beta_0 +  2 s_n M_n' \beta_n
  e^{-n\lambda_k + n\epsi} \beta_0)$, which is $\D^k(\tau_n)$. That
  proves $S_{\hat x}(P_n) \subset \D^k(\tau_n)$. \fin

Now let us restate propositions A and B of section \ref{pmk}. We parametrize the family  $\LL_n$ of polydiscs by $(i,\al) \in
\{ 1,\ldots,k \} \times \D(\tau_n)$. More precisely, let $\Ld_n^{i,\al} : \D^{k-1} \to
\D^k(\tau_n)$ be defined as $\Ld_n^{i,\al}(v_1,\ldots,v_{k-1}) = (v_1 \tau_n , \ldots , \al ,
\ldots , v_{k-1}
\tau_n )$, 
where $\al$ stands at the $i$-th coordinate. Pulling back  $\Ld_n^{i,\al}$ by $S_{\hat x}$, we set $L_n^{i,\al} := S_{\hat
  x}^{-1} \circ \Ld_n^{i,\al}$. By lemma \ref{nouv}(3), that polydisc
satisfies  $L_n^{i,\al} : \D^{k-1} \to B_x(\tau_n)$. \\

Proposition A now take the following form.

\begin{proposition}{A:}
For every $(i,\al)  \in  \{ 1,\ldots,k \} \times \D(\tau_n)$, $\Vol f^n \circ L_{n}^{i,\al}  \leq  d^{(k-1)(n-q_n)}$.  
\end{proposition}

\noindent Before dealing with proposition B, let us introduce the collection    
\begin{equation*} 
    \PP'_n(x) := \left \{ \, f^{-n}_{\hat y_n}( B_p(s_n/2) ) \, , \,  y \in B_x^{\Om_\epsi}(\rho_n) \, , \,  p \in \EE_{s_n}(y_n)  \right \} .
\end{equation*} 
It satisfies $\cdb \PP'_n(x) = \cdb \PP_n(x)$ and its sets are pairwise
disjoint. Given $P_n \in \PP'_n(x)$, for simplicity we denote $\Vol f^n  \left(
L_n^{i,\al} \cap P_n \right)$ for the volume of $f^n \circ
L_n^{i,\al}$ restricted to $(L_n^{i,\al})^{-1}(P_n)$. Observe that it has
multiplicity 1 since $f^n$ is injective on $P_n$. Proposition B is
restated as follows.

\begin{proposition}{B:}
There exists a subset $\Lambda_n \subset \D(\tau_n)$ satisfying: $\cdb \Lambda_n \leq e^{20 n \epsi}$ and for every $P_n \in \PP'_n(x)$, there is $(i,\al)(P_n) \in \{
 1,\ldots,k \} \times \Lambda_n$ such that: 
\begin{equation*}\label{rbn}
 \Vol f^n  \left( L_n^{(i,\al)(P_n)} \cap P_n \right)  \geq (s_n)^{k-1} . 
\end{equation*} 
\end{proposition}

\noindent Let us see how we deduce (\ref{tha}), thus completing the
proof of theorem B. Since the sets of $\PP'_n(x)$ are pairwise disjoint, we have:
\[ \sum_{P_n \in \PP'_n(x)} \Vol f^n  \left( L_n^{(i,\al)(P_n)} \cap P_n
\right)  \leq   \sum_{i=1}^k \sum_{\al \in \Lambda_n} \Vol f^n  \circ L_n^{i,\al} .   \] 
That implies $\cdb \PP_n(x) \cdot  (s_n)^{k-1}  \leq k \,
\cdb \Lambda_n \cdot d^{(k-1)(n-q_n)}$. Then (\ref{tha}) follows using $s_n = r_0 e^{-8 n\epsi}$ and $\cdb \Lambda_n
\leq e^{20 n \epsi}$.

\section{Proof of propositions A and B (stated in §\ref{cedec})}\label{proay}

\subsection{Proof of proposition A}\label{proa}

We denote by $\bar L_n^{i,\al}$ the
extension of $L_n^{i,\al}$ to the polydisc $\D^{k-1}(2)$, it satisfies
$\bar L_n^{i,\al} \subset B_x(2 \tau_n)$. We
set $\si_{q_n} := f^{q_n} \circ L_n^{i,\al}$ and $\bar \si_{q_n} :=
f^{q_n} \circ \bar L_n^{i,\al}$. According to section \ref{pmk}, proposition A is a consequence of the
growth lemma combined with the following lemma:

\begin{lem}\label{paire}
The polydisc $\bar \si_{q_n}$ is bounded.
\end{lem}

\proof  We have to show that  $\bar \si_{q_n}$ is included in some
 $B_j$. With no loss of generality, we can assume that every
 ball of radius $r_0$ in $\Pj^k$ is contained in some $B_j$. For simplicity we denote $q := q_n$. Observe that it suffices to prove 
\begin{equation}\label{compp} 
  \D^k ( 2 \beta_0 \tau_n )  \subset  R^{q}_{\hat x_{q}}  \circ S_{\hat
 x_{q}} ( B_{x_{q}} ( r_{q} ) ).
\end{equation}
Indeed, that inclusion implies using $R^{q}_{\hat x_{q}}  \circ S_{\hat
 x_{q}} = S_{\hat x}
 \circ f_{\hat x_{q}}^{-q}$ and $r_{q} \leq r_0$:
\[  f^{q} \circ S_{\hat x}^{-1} ( \D^k ( 2 \beta_0 \tau_n ) )  \subset B_{x_{q}} ( r_0 ). \]  
The conclusion then follows from (see lemma \ref{nouv}(3)
for the last inclusion): 
\[ \bar \si_{q} =
f^{q} \circ \bar L_n^{i,\al} \subset f^{q} ( B_x(2 \tau_n) )
\subset  f^{q} \circ S_{\hat x}^{-1} ( \D^k ( 2 \beta_0 \tau_n )). \]
Thus it remains to show (\ref{compp}). Lemma \ref{nouv}(4,3) yields:
\[ \D^k ( {M_{q}'}^{-1} e^{- q\lambda_1 - q \epsi} \cdot r_{q}  ) \subset  R^{q}_{\hat x_{q}}  ( \D^k ( r_{q} )) \subset     R^{q}_{\hat x_{q}}  ( S_{\hat x_{q}} ( B_{x_{q}} ( r_{q} ) )  ) . \] 
Using $q \lambda_1 \leq n \lambda_k$ (which implies $q \leq n$), we obtain: 
\[   {M_{q}'}^{-1} e^{- q\lambda_1 - q \epsi} \cdot r_{q}  =  r_0
{M_0'}^{-1} e^{- q\lambda_1 - 3 q \epsi} \geq  e^{- n\lambda_k - 4 n
  \epsi}  \geq 2\beta_0^2 ( 1 + 2  \beta_0 M_0'
) r_0 e^{-n\lambda_k -5n\epsi}, \] 
which is equal to $2 \beta_0 \tau_n$. \finsec

\subsection{Proof of proposition B}\label{proofb}

We set $\eta_n := s_n  e^{-n\lambda_k - 4n\epsi} / 4$ and define
$\Lambda_n$ as a maximal $\eta_n$-separated set in $\D(\tau_n e^{n \epsi})$. 
We have $\cdb \Lambda_n \leq (\tau_n
e^{n \epsi})^2 / \eta_n^2 \leq 
e^{20n\epsi}$. Let us fix $P_n  \in \PP'_n(x)$ for the remainder of the
section and show:
\begin{equation} \label{zzp}
  \exists \, (i,\al)(P_n) \in \{1,\ldots,k \}
 \times \Lambda_n \ , \  \Vol f^n  \left( L_n^{(i,\al)( P_n)} \cap P_n \right)  \geq (s_n)^{k-1} . 
\end{equation}
Let also $y \in B_x^{\Om_\epsi}(\rho_n)$ and $p \in
\EE_{s_n}(y)$ such that $P_n = f^{-n}_{\hat y_n} (P) = f^{-n}_{\hat y_n}( B_p(s_n/2))$.

\subsubsection{Definition of $(i,\al)(P_n)$}\label{opo}

We define $1 \leq i(P_n) \leq k$ to be the unique element satisfying $y \in
Y_x^{i(P_n)}$ (see subsection \ref{defam}). For simplicity we denote $j := i(P_n)$. We now define $\al(P_n) \in \Lambda_n$. Since
$S_{\hat x}(P_n) = S_{\hat x} \circ f^{-n}_{\hat y_n} (P)  \subset
\D^k(\tau_n)$ (lemma \ref{where}), then $\pd_n := S_{\hat x} \circ f^{-n}_{\hat
  y_n}(p)$ lies in $\D^k(\tau_n)$. In particular we have $\pi_j({\pd}_n) \in
\D(\tau_n)$ and $\D_{\pi_j({\pd}_n)}(\eta_n)
\subset \D(\tau_n e^{n \epsi})$. In order to find some $\al(P_n) \in
\Lambda_n$ satisfying (\ref{zzp}), we shall prove:
\begin{equation}\label{ff}
 \forall \al \in \D_{\pi_j({\pd}_n)}(\eta_n) \ , \ \Vol f^n \left(L_{n}^{j,\al} \cap P_n  \right) \geq  (s_n)^{k-1}. 
\end{equation}
Then we take for $\al(P_n)$ any element in $\Lambda_n \cap
\D_{\pi_j({\pd}_n)}(\eta_n)$: that set is not empty since $\Lambda_n$
is a maximal $\eta_n$-separated set in $\D(\tau_n e^{n \epsi})$. That
shows theorem B. \\

We deduce (\ref{ff}) from the following claim. It relies on a precise geometrical description of the inverse branches, due to the normalization theorem. We set $Q := B_p(s_n
/4)$, $Q_n := f^{-n}_{\hat y_n}(Q)$ and identify the polydisc $L_{n}^{j,\al}$ with its source $\D^{k-1}$.

\begin{claim}{:}
For every $\al \in \D_{\pi_j({\pd}_n)}(\eta_n)$, 
\begin{enumerate}
\item[(a)] $L_{n}^{j,\al}$ intersects $Q_n$,
\item[(b)] the slice $P_n \cap L_{n}^{j,\al}$ is a domain in $\D^{k-1}$ with boundary in $\partial P_n$. 
\end{enumerate}
\end{claim}

Let us see how we infer (\ref{ff}). Let $a \in Q_n \cap L_{n}^{j,\al}$. Since $f^n(a) \in Q$, we have $Q' := B_{f^n(a)}( s_n /4) \subset P = B_p(s_n /2)$.
Hence $\Sigma := f^n(P_n \cap
L_{n}^{j,\al})$ satisfies $\Sigma  \subset P$ and
$\partial \Sigma \subset\partial P$ (the map $f^n : P_n \to P$ 
is a biholomorphism). Therefore $\Sigma  \cap  Q'$ is an immersed polydisc containing $f^n(a)$ (the center of $Q'$) with boundary in $\partial Q'$. The Lelong inequality \cite{L} then implies $\Vol \left(
\Sigma  \cap  Q' \right) \geq (s_n)^{k-1}$ up to a multiplicative
constant. That gives (\ref{ff}) and completes the proof of
theorem B.

\subsubsection{Proof of the claim}\label{popl}

Let us denote $\psi := \psi_{x,y}$. For every $s \leq s_n$, we set $\eta := s e^{-n\lambda_k - 4n\epsi}$, $\Ad := S_{\hat y_n}( B_p(s))$ and
$\Ad_n := \psi_{x,y}
\circ R^n_{\hat y_n} (\Ad)$. For simplicity we assume that
  $\Ad = \D^k_{\pd} (s)$, where $\pd :=
S_{\hat y_n}(p)$ (see lemma \ref{nouv}(3)).
   For any $\tilde u = (u_1,\ldots,u_{k-1}) \in \D^{k-1}$, we define
  $\vd_{\tilde u} : \D \to \Ad$ by $\vd_{\tilde u}(t) := \pd  + s
  ({\tilde u} , t)$. The claim is a consequence of the next proposition applied with $s =
s_n / 2$ (for the item (b)) and $s_n / 4$ (for the item (a)). 
\begin{prop}\label{jhn}  For every ${\tilde u} \in \D^{k-1}$,
\begin{enumerate}
\item  $\psi_{x,y} \circ
R^n_{\hat y_n} (\vd_{\tilde u})$ is a graph over the
  $j$-axis,  
\item its $\pi_j$-projection $\wwd_{\tilde u}^n := \pi_j \circ
\psi_{x,y}  \circ R^n_{\hat y_n} (\vd_{\tilde u})$ contains the disc 
$\D_{\pi_j({\pd}_n)}(\eta)$.
\end{enumerate}
\end{prop}

We need the following lemma for proving proposition \ref{jhn}.

\begin{lem}\label{resn} For every ${\tilde u} \in \D^{k-1}$,
  $\wwd_{\tilde u}^n : \D \to \C$ satisfies :
\begin{enumerate}
\item $\wwd_{\tilde u}^n(0) \in \D_{\pi_j(\pd_n)} ( s   e^{-n\lambda_{k-1}+3n\epsi} )$,
\item $\forall t \in \D$, $\abs{ {\wwd_{\tilde u}^n}'(t) - {\wwd_{\tilde u}^n}'(0)  }  \leq
  s   e^{- 2n\lambda_k + 3 n \epsi}$,
\item $\abs{ {\wwd_{\tilde u}^n}'(0) } \geq s    e^{-n\lambda_k-2n\epsi}$. 
\end{enumerate} 
\end{lem}

\proo \textsc{of proposition \ref{jhn}:} For the point 1, it suffices to verify that $\wwd^n_{\tilde u} = \pi_j \circ
\psi_{x,y}  \circ R^n_{\hat y_n} (\vd_{\tilde u})$ is
injective. Let $\varphi :=  ( \wwd_{\tilde u}^n - \wwd_{\tilde u}^n(0)
) / {\wwd_{\tilde u}^n}'(0) - \Id$. We get from lemma \ref{resn}(2,3):
\[ \forall t \in \D \ , \ \abs{ \varphi'(t) }  = { \abs{ {\wwd_{\tilde u}^n}'(t) - {\wwd_{\tilde u}^n}'(0)  }  \over \abs{ {\wwd_{\tilde u}^n}'(0)  } } \leq  { e^{-2n \lambda_k + 3n \epsi}
  \over   e^{-n\lambda_k - 2n\epsi} }  = e^{-n\lambda_k  + 5 n \epsi} .  \]
This implies $\Lip (\varphi) \leq 1/2$, hence $\Id +
\varphi$ and $\wwd_{\tilde u}^n$ are injective on $\D$. Let us prove
the point 2. Since $\Lip (\varphi) \leq 1/2$ and
$\varphi(0)=0$, we have $\abs{(\Id + \varphi)(t)} \geq \abs{t} -
\abs{\varphi(t)} \geq \abs{t} / 2$. That yields $\abs{\wwd_{\tilde u}^n(t) - \wwd_{\tilde u}^n(0)} \geq
\abs{{\wwd_{\tilde u}^n}'(0)} / 2$ for every $t \in \Sg^1$. Then lemma
\ref{resn}(3) implies:
\[ \forall t \in \Sg^1 \ , \  \abs{ \wwd_{\tilde u}^n(t) - \wwd_{\tilde u}^n(0)} \geq  s    e^{-n\lambda_k-3n\epsi},  \]
which yields $\D_{\wwd_{\tilde u}^n(0)}
\left(s   e^{-n\lambda_k -3n\epsi} \right ) \subset
\wwd_{\tilde u}^n$ by Jordan's theorem. We deduce from lemma \ref{resn}(1) that:  
\begin{equation*}
  \forall \tilde u \in \D^{k-1} \ , \ \D_{\pi_j({\pd}_n)} \left( s   e^{-n\lambda_k - 3n\epsi} - s    e^{-n\lambda_{k-1} + 3n\epsi} \right) \subset  \D_{\wwd_{\tilde u}^n(0)} \left(s    e^{-n\lambda_k-3n\epsi} \right) .
\end{equation*}
Finally, the left hand side contains $\D_{\pi_j({\pd}_n)} (\eta) = \D_{\pi_j({\pd}_n)}
 (s e^{-n \lambda_k - 4n\epsi})$.  \finsec

\subsubsection{Proof of lemma  \ref{resn}}\label{rtt}

We shall use the algebraic properties of resonant maps (namely lemma
\ref{bnormb}). For every $(\tilde u,t) \in \D^{k-1} \times \D$, we
denote $\Ad(\tilde u,t) := \vd_{\tilde
  u}(t) =  \pd + s (\tilde
u,t)$ and $z := \Ad(\tilde u,t)$. We also denote:
\[ \vd_{\tilde
 u}^n(t) := R^n_{\hat y_n} \circ  \vd_{\tilde u}(t)  \ \ \textrm{ and }  \ \ \hd^n (\tilde u) :=  R^n_{\hat y_n} \circ \Ad(\tilde
 u, 0 ) . \] 
We have therefore $\pd_n = \psi \circ R^n_{\hat y_n}(\pd) = \psi
\circ \hd^n (0)$. Observe that $\Ad \subset
\D^k$ (lemma \ref{gds}(2)) implies $\vd_{\tilde
 u}^n \subset \D^k(M_n' e^{-n\lambda_k + n\epsi}) \subset \D^k(R)$ (lemma
 \ref{nouv}(4)). One also obtains from
 the very definition of resonant maps (see (\ref{ztop}), subsection \ref{kop}):
\begin{equation}\label{topol} 
{\vd_{\tilde u}^n}'\, (t) = s  \cdot {\partial  R^n_{\hat y_n} \over \partial z_k} (z) \ \ , \ \ {\vd_{\tilde
 u}^n}''\, (t) = s^2 \cdot {\partial^2  R^n_{\hat y_n} \over \partial
 z_k^2} (z) \ \ \textrm{ and}  \ \ \pi_k \circ \hd^n \equiv  a_{k,n}(\hat y_n) \cdot \pi_k(\pd).
\end{equation} 
We deduce from the last observation:
\begin{equation}\label{top}
\norm{ d_{\tilde u}
\hd^n } = \norm{ \tilde \pi \circ  d_{\tilde u} \hd^n } = \norm{ \tilde \pi
   \circ  d_z R^n_{\hat y_n} \circ  d_{(\tilde u,0)} \Ad  }   = s \norm{ \tilde \pi \circ d_z
 R^n_{\hat y_n} }. 
\end{equation}
Finally let us recall that $\wwd_{\tilde u}^n  = \pi_j \circ \psi \circ
\vd_{\tilde u}^n$.\\

1 - $\wwd_{\tilde u}^n(0) \in \D_{\pi_j(\pd_n)} ( s   e^{-n\lambda_{k-1}+3n\epsi} )$.

\vskip 3 pt

\noindent We have $\wwd_{\tilde u}^n(0) \in \pi_j \circ \psi \circ
\hd^n(\D^{k-1})$ and $\pi_j(\pd_n) = \pi_j \circ \psi \circ
\hd^n(0)$. Moreover (\ref{top}) yields for every $\tilde u \in \D^{k-1}$:
\[ \norm{ d_{\tilde u}(\pi_j \circ \psi  \circ {\hd^n})  } \leq \norm {
  d_{\hd^n(\tilde u)} \psi } \, \norm{ d_{\tilde u} \hd^n } = s  \norm {
  d_{\hd^n(\tilde u)} \psi } \, \norm{ \tilde \pi \circ d_z
 R^n_{\hat y_n} } \]
which is less than $s \beta_0 M_n' e^{-n\lambda_{k-1} + n \epsi} \leq
  s   e^{-n\lambda_{k-1} + 3 n \epsi}$ (lemmas \ref{bnormb}(2)
  and \ref{controun}(2)). That proves the point 1.\\

2 - $\forall t \in \D$, $\abs{ {\wwd_{\tilde u}^n}'(t) - {\wwd_{\tilde u}^n}'(0)  }  \leq
  s   e^{- 2n\lambda_k + 3 n \epsi}$.

\vskip 3 pt

\noindent Since $\wwd_{\tilde u}^n  = \pi_j \circ \psi \circ
\vd_{\tilde u}^n$, it suffices to verify that $\phi_{\tilde u}^n := ( \psi \circ \vd_{\tilde u}^n ) ' -   ( \psi \circ \vd_{\tilde u}^n )'(0)$ satisfies $\abs {\phi_{\tilde u}^n} \leq s  e^{-2n\lambda_k +3n \epsi}$. Let us write for every $t \in \D$ :  
\[ \phi_{\tilde u}^n(t)   =   \left( d_{{\vd_{\tilde u}^n}(t)} \psi  - d_{{\vd_{\tilde u}^n}(0)} \psi \right) ({\vd_{\tilde u}^n}'(t)) +  \left( d_{{\vd_{\tilde u}^n}(0)} \psi \right) \left( {\vd_{\tilde u}^n}'(t) - {\vd_{\tilde u}^n}'(0) \right). \]
Using lemma \ref{controun}(2,3), we obtain for every $t \in \D$:
\begin{equation*}\label{jjk}
  \abs{ \phi_{\tilde u}^n(t) }  \leq  \ga  \abs{ {\vd_{\tilde u}^n}(t) - {\vd_{\tilde u}^n}(0)} \abs{{\vd_{\tilde u}^n}'(t)} + \beta_0   \abs{  {\vd_{\tilde u}^n}'(t) - {\vd_{\tilde u}^n}'(0) } \leq  \ga  \abs{  {\vd_{\tilde u}^n}'
}_{\infty,\D}^2 +  \beta_0  \abs{  {\vd_{\tilde u}^n} '' }_{\infty,\D}.
\end{equation*}
We deduce using (\ref{topol}) and lemma \ref{bnormb}(3,4):
\[  \abs{ \phi_{\tilde u}^n(t) }  \leq
\ga s^2  \max \{ M_n' e^{-n\lambda_{k-1} + n\epsi}   \, , \,  e^{-n \lambda_k + n\epsi} \} ^2 + \beta_0 s^2 M_n' e^{- 2 n 
  \lambda_k +  n \epsi}  \leq s  e^{-2n\lambda_k + 3 n \epsi} . \]
That proves the point 2.\\

3 - $\abs{ {\wwd_{\tilde u}^n}'(0) } =  \abs{ \left( \pi_j
\circ d_{{\vd_{\tilde u}^n}(0)} \psi \right) ({\vd_{\tilde
    u}^n}'(0)) } \geq s  
  e^{-n\lambda_k-2n\epsi}$.

\vskip 3 pt

\noindent The line (\ref{topol}) and lemma \ref{bnormb}(2,3) yield  for ${\vd_{\tilde
    u}^n}'(0) \in \C^k$ :
\begin{equation*}\label{un}
 \abs { \tilde \pi ( {\vd_{\tilde u}^n}'(0) ) }  \leq s   M_n' e^{-n
  \lambda_{k-1} + n\epsi} \ \ \textrm{
  and } \ \  \abs{ \pi_k ( {\vd_{\tilde u}^n}')(0)}
  \geq s    e^{-n \lambda_k - n\epsi}. 
\end{equation*}
Now lemmas \ref{controun}(2) and \ref{ggp} imply (use $y \in Y^j_x$ for the second inequality):
\begin{equation*}\label{deux}
  \forall 1 \leq i \leq k - 1  \ , \ \abs{ (\pi_j \circ d_{\vd_{\tilde u}^n(0)} \psi) (c_i)  } \leq
  \beta_0 \ \ \textrm{
  and } \ \ \abs{ ( \pi_j \circ  d_{\vd_{\tilde u}^n(0)} \psi )
  (c_k)} \geq 1/(2\beta_0). 
\end{equation*}
We deduce $\abs{  {\wwd_{\tilde u}^n}'(0) } \geq  s \left( (2\beta_0)^{-1}
  e^{-n\lambda_k - n\epsi} - \beta_0 M_n'  e^{-n\lambda_{k-1} +n \epsi}
  \right) \geq s  e^{-n\lambda_k-2n\epsi}$, completing the proof of
  lemma \ref{resn}.

\section{Appendix}

Let $f$ be a holomorphic endomorphism of $\Pj^k$ with degree $d \geq
2$. Let $\om$ be the Fubini-Study $(1,1)$ form on $\Pj^k$. For every
holomorphic polydisc $\eta : \D^l (r) \to \Pj^k$, we define $\Vol f^m
\circ \eta := \int_{\D^l(r) } \eta^* {f^m}^* \om^l$. We recall that $\{ B_j \, , \, j \in J
\}$ is a finite covering of $\Pj^k$ which consists of open sets bounded in the affine charts. We say that $\eta$ is bounded if
the image of that polydisc is contained in some $B_j$. 
\begin{fact}{:}\label{aire}
Let $1 \leq l \leq k$ and $\eta : \D^l (2) \to \Pj^k$ be a bounded polydisc. Then \[ \forall m \geq 1 \, , \, \Vol  f^{m} \circ
\eta_{\vert \D^l}  \leq d^{l m}. \] 
\end{fact}
The proof relies on the Green current of $f$, which is the closed positive
$(1,1)$ current on $\Pj^k$ defined by $T = \lim_{n \to \infty} {1 \over d^n} {f^n}^* \om$.
 That current satisfies $f^* T = d T$ and $T = \om - dd^c \varphi$ for some continuous function $\varphi
 : \Pj^k \to \R$. Iterating that identity, we obtain:
\begin{equation}\label{popol}
\forall m \geq 1 \  ,  \  {f^m}^* \om  = d^m T + dd^c (\varphi \circ f^m). 
\end{equation}
We refer to the article of Dinh-Sibony (\cite{DS}, section 1.2) for
more details about the Green current. In order to prove the lemma, we
shall use an induction concerning the mass of $T^i \wedge {f^m}^*
\om^j$. Note that a similar induction was employed by Dinh to estimate
the local entropy outside the support of the current $T^i$ (see \cite{D},
theorem 2.1). In the sequel $\om_0$ stands for the $(1,1)$ form $dd^c
\abs{\,.\,}^2$ which induces the standard metric  on $\C^k$.

\proo \textsc{of the growth lemma:} It follows from Cauchy's estimates that the family of bounded polydiscs $\D^l (2) \to \Pj^k$ has bounded derivatives on $\D^l (3/2)$, say by $1$. We deduce that for any such polydisc $\eta$ and any positive current $S$ on $\Pj^k$ of bidegree $(s,s)$ (with $s \leq l$): 
\begin{equation}\label{otra}
  \forall \rho \leq 3/2  \ , \  0 \leq \int_{\D^l(\rho)} \eta^* S \wedge \eta^* \om^{l-s} \leq  \int_{\D^l(\rho)} \eta^* S \wedge \om_0^{l-s}.
\end{equation}
Let us fix $\eta : \D^l (2) \to \Pj^k$ and denote by $\norm{S}_\rho := \int_{\D^l(\rho)} \eta^* S \wedge \om_0^{l-s}$. We shall prove for any $1 \leq q \leq l$ and $0 \leq r \leq q$:
\[  (H_{q,r}) \ : \ \exists c_{q,r} \geq 1 \ , \  \exists \rho_{q,r} \in ] 1,3/2 [ \ , \ \forall m \geq 0 \ , \ \norm{T^{q-r} \wedge {f^m}^* \om^r  }_{ \rho_{q,r}} \leq  c_{q,r} \, d^{mr} . \]
The lemma then follows by taking $S = {f^m}^*\om^l$ and $s=l$, and by using (\ref{otra}) and $(H_{l,l})$. \\

Let us establish $(H_q):=$ ``$(H_{q,r})$ holds for any $0 \leq r \leq q$'' by induction on $q$. Observe that $(H_{q,0})$ obviously holds for any $1 \leq q \leq l$. Hence it suffices to verify $(H_{1,1})$ to end the proof of $(H_1)$.  Let $1 < \rho_{1,1} < \tau_{1,1} < \rho_{1,0} < 2$ and $\chi$ be a cut-off function with support in $\D^l ( \tau_{1,1} )$ such that $\chi \equiv 1$ on $\D^l (\rho_{1,1})$. We deduce from (\ref{popol}):  
\begin{equation*}\label{mpp} 
 \norm{  {f^m}^* \om }_{\rho_{1,1}} = d^m  \norm{T}_{\rho_{1,1}}  + \int_{\D^l(\rho_{1,1})}  dd^c (\varphi \circ f^m \circ \eta) \wedge \om_0^{l-1} =: d^m \, \norm{T}_{\rho_{1,1}} + A_m . 
\end{equation*}
On one hand $\norm{T}_{\rho_{1,1}} \leq  c_{1,0}$ from $(H_{0,0})$. On the other hand Stokes' theorem implies up to some multiplicative constant:
\[ A_m \leq \int_{\D^l(2)} \chi \cdot dd^c  (\varphi \circ f^m \circ \eta) \wedge \om_0^{l-1} = \int_{\D^l(2)} \varphi \circ f^m \circ \eta \cdot dd^c \chi \wedge \om_0^{l-1} \leq  \norm {\varphi}_\infty \norm{\chi}_{C^2}.  \]
Hence there exists $c_{1,1} \geq 1$ such that $\norm{  {f^m}^* \om }_{\rho_{1,1}} \leq
c_{1,1} d^m$, which proves $(H_1)$. \\

Assume now that $(H_q)$ holds
for $1 \leq q \leq l-1$, and let us prove $(H_{q+1})$. For that purpose, we show
$(H_{q+1,r})$ by induction on $r$. Given $0 \leq r \leq q$, we shall deduce
$(H_{q+1,r+1})$ from $(H_{q,r})$ and $(H_{q+1,r})$. Let us set $1 <
\rho_{q+1,r+1} < \tau_{q+1,r+1} < \min \{ \rho_{q,r} , \rho_{q+1,r}
\}$ and let $\chi$ be a cut-off function with support in $\D^l ( \tau_{q+1,r+1} )$ such that $\chi \equiv 1$ on $\D^l (\rho_{q+1,r+1})$. We obtain using (\ref{popol}):
\begin{equation}\label{mppp} 
  T^{q+1-(r+1)} \wedge {f^m}^* \om^{r+1}    =  T^{q-r} \wedge {f^m}^* \om^r \wedge \big ( d^m T + dd^c (\varphi \circ f^m)  \big) = d^m \, S_1 + S_2 ,                      
\end{equation}    
where $S_1 := T^{q+1-r} \wedge {f^m}^* \om^r$ and $S_2 := T^{q-r}
\wedge {f^m}^* \om^r \wedge dd^c (\varphi \circ f^m)$. Now
$(H_{q+1,r})$ and $(H_{q,r})$ respectively imply (use  Stokes' theorem as before
for the second line): 
\[ d^m \, \norm{S_1}_{\rho_{q+1,r+1}} \leq  d^m \, \norm{S_1}_{\rho_{q+1,r}} \leq  c_{q+1,r} \, d^{m(r+1)}  , \]
 \[ \norm{S_2}_{\rho_{q+1,r+1}} \leq \norm {\varphi}_\infty \norm{\chi}_{C^2} \norm { T^{q-r} \wedge {f^m}^* \om^r }_{\rho_{q,r}} \leq  \norm {\varphi}_\infty \norm{\chi}_{C^2}  c_{q,r} \, d^{mr} . \]
Using (\ref{mppp}) we get $\norm{ T^{q+1-(r+1)} \wedge {f^m}^*
  \om^{r+1}  }_{\rho_{q+1,r+1}} \leq c_{q+1,r+1} \, d^{m(r+1)}$ for
some $c_{q+1,r+1} \geq 1$. That completes the proof of the growth lemma. \finsec

\vspace{1 cm}

{\footnotesize C. Dupont}\\
{\footnotesize Universit\'e Paris XI - Orsay}\\
{\footnotesize CNRS UMR 8628}\\
{\footnotesize Math\'ematiques, B\^at. 425}\\
{\footnotesize F-91405 Orsay Cedex, France}\\
{\footnotesize christophe.dupont@math.u-psud.fr}

\end{document}